\documentclass[11pt]{article}
\usepackage[english]{babel}
\usepackage{cite}
\usepackage{amsmath,amssymb}
\usepackage{epsfig}

\begin{document}
\thispagestyle{empty}
\centerline{DESY 17-219 \hfill ISSN 0418-9833}
\centerline{May 2018 \hfill}

\vspace*{2.0cm}

\begin{center}
\boldmath
 {\large \bf
  q-derivatives  of  multivariable q-hypergeometric function  with respect to their parameters
 }
\unboldmath
\end{center}
 \vspace*{0.8cm}

\begin{center}
{\sc Vladimir~V.~Bytev,$^{a}$\footnote{E-mail: bvv@jinr.ru}
Pengming Zhang ,$^{b,}$\footnote{E-mail: zhangpm5@mail.sysu.edu.cn}
} \\
 \vspace*{1.0cm}
{\normalsize $^{a}$ Joint Institute for Nuclear Research,} \\
{\normalsize $141980$ Dubna (Moscow Region), Russia}\\
\bigskip
{\normalsize $^{b}$  School of Physics and Astronomy, Sun Yat-sen University, } \\
{\normalsize  519082 Zhuhai, China,}
\end{center}

\begin{abstract}
We consider the q-derivatives   of the Srivastava and Daoust basic multivariable hypergeometric function with respect to
the parameters. This function embodies a entire number  of various q-hypergeometric series of one and several variables.
 Explicit equations are given for general case of summation indexes with positive real  coefficients. As an example
derivatives of q-analog of non-confluent Horn type hypergeometric function  $H_3$ is presented.
\end{abstract}

\newpage

\section{Introduction}
\label{Intro}

Basic hypergeometric functions, namely q-analogue generalizations of ordinary hypergeometric series, has a long  history \cite{book1}.
Contrary to the generalized hypergeometric series, where  the ratio of  successive terms  is a rational function of n,  the ratio of successive terms in basic hypergeometric series is a rational function of $q^n$.

The basic hypergeometric series $ _2\phi_1(q\alpha,q\beta;q\gamma;q,x)$, q-extension of the Gauss hypergeometric function $_2F_1$ was first considered by Eduard Heine in 1846.
The last could be restored from q-analogue by applying parameters equal $q^{\alpha_i}$ and considering  the limit when the base $q$ tends to unity.

Renewed interest in hypergeometric-type series containing digamma functions has
emerged in connection with derivatives of hypergeometric functions with respect
to their parameters.
The first derivatives for some special values of the parameters were already
known a long time ago \cite{firstDer1,firstDer2,firstDer3}.
Later on, Ancarani et al.\ have found in a series of papers
\cite{Ancarani0,Ancarani1,Bujar} the derivatives of Gaussian hypergeometric
functions and some derivatives of two-variable series, namely the Appell
series and four degenerate confluent series~\cite{Ancarani2}.
Moreover, it has been shown that the first derivatives of generalized
hypergeometric functions are expressible in terms of Kamp{\'e} de F{\'e}riet
functions, and, with the same technique, derivatives of the Appell
hypergeometric function have been obtained in Ref.~\cite{Sahai}.

In the paper \cite{Bytev:2017jmx} the general case of hypergeometric function derivative with respect to the parameters was considered.
Namely, the case of multi-summation index with arbitrary natural coefficients was calculated. It was shown, that the derivatives with respect 
to the upper and lower parameters of the (Srivastava-Daoust) generalized Lauricella series are expressed as a final sum of the  generalized Lauricella series
with shifted parameters.

In quantum field theory, one has to calculate higher-order Feynman diagrams
for quantum corrections to scattering processes.
For the evaluation of the  dimensionally regularized Feynman
integrals in $D=4-2\varepsilon$ dimensional space time, one has to construct
the expansion in $\varepsilon$. In the   \cite{Bytev:2017jmx} it was concluded that the $\varepsilon$ expansions
of Feynman integrals at any order are expressible in terms of Horn-type
hypergeometric functions: the $n$-th term of the $\varepsilon$ series can be expressed as a
Horn-type hypergeometric function in $n+m$ variables,
where $m$ is the number of summations in the Horn-type representation of the
Feynman integral.
The region of convergence of any of these parameter derivatives, i.e., the
coefficients in the $\varepsilon$ expansion, and the initial Feynman integral
are the same.

Classes of the Feynman integrals that could be expressed in the terms of multiple polylogaritms (MP) \cite{levin} are
mainly related to the massless cases. 
More interesting massive case, for example  massive  sunrise and kite integrals could be expressed through elliptic generalization of MP (EP), see 
\cite{Adams:2015ydq},\cite{Adams:2016xah}. At \cite{Passarino:2016zcd} it was shown, that EP can be written as a multiple integral  of basic hypergeometric function convolution. The utilization of q-difference equations and q-contiguous relations gives one the possibility to establish a connection with basic hypergeometric functions and find  efficient algorithm  for analytical continuation  and numerical evaluation of EP.  
By analogy with paper \cite{Bytev:2017jmx}  the knowledge of  of q-basic function derivatives over parameters can help us to find out the properties of   Feynman integral expansion over dimensional parameter and make stable numerical estimations.

The generalization of hypergeometric series derivatives, the q-derivatives of basic hypergeometric series with respect to the parameters was considered in the papers \cite{qref1,qref2,qref3,qref4,qref5}.  The case of q-hypergeometric series $_r\Phi_s$ (q-extension of generalized hypergeometric series) was studied at \cite{qref1},
later the derivatives of Appell functions with respect to parameters and q-derivatives of 3-variable q-Lauricella functions \cite{qref2} as well as k-
variable q-Lauricella functions and q-Kampe de Feriet function \cite{qref3,qref4} was calculated. Some special cases of q-derivative  of Srivastava's general triple q-hypergeometric series
with respect to its parameters $H_{A,q},H_{B,q},H_{C,q}$ are considered in  \cite{qref5}.

These multiple basic hypergeometric series  considered above are special cases of the q-extension of generalized Lauricella series in $n$ variables (\ref{Eq:Sriv:Def}).  At this paper we would like to present  the q-derivatives of the basic (Srivastava-Daoust) generalized Lauricella series with respect to the upper and lower parameters and generalize results obtained 
at \cite{qref1,qref2,qref3,qref4,qref5}. Here we consider the case when coefficients  $\theta_1^{{(1)}},\psi_1^{{(1)}},\phi_1^{{(1)}},\delta_1^{{(1)}},\dots, \theta_A^{{(n)}},\psi_C^{{(n)}},\phi_{B^{(n)}}^{{(n)}},\delta_{D^{(n)}}^{{(n)}}$ of multi summation indexes in (\ref{eq:def:omeg}) are positive reals, contrary to the considered special cases, where these coefficient are implied equal zero or unity.

This paper is organized as follows.
We begin in Sec.~\ref{Sect:Definitions} by considering definitions for q-extension of generalized Lauricella hypergeometric function and some basic definitions of the q-special function theory, namely Jackson derivative, q-Pochhammer symbol and   q-bracket.

Next, Sec.~\ref{sec3} is devoted to the q-derivatives in the case of double
summation index parameters with real positive coefficients, namely as an example the  q-analog of non-confluent Horn type hypergeometric function $H_3(a,b,c,x,y)$
q-derivatives with respect to the upper and lower parameters are considered.

Our main result are presented in  Sec.~\ref{sec:Main}, where we consider the  q-derivatives  of  q-extension of the (Srivastava-Daoust) generalized Lauricella series in $n$ variables
with respect to the upper and lower parameters. The cases of multiple summation indexes with positive real coefficient are considered.

\section{Definitions}
\label{Sect:Definitions}

Q-extension of the (Srivastava-Daoust) generalized Lauricella series in $n$ variables, given by Srivastava \cite{Sriv1}
is defined:
\begin{eqnarray}
\lefteqn{F^{A:B^{(1)};\dots; B^{(n)}}_{C:D^{(1)};\dots; D^{(n)}} \left( \begin{array}{c}
x_1 \\
\vdots\\
x_n
\end{array} \right)}
\nonumber\\
&=&F^{A:B^{(1)};\dots; B^{(n)}}_{C:D^{(1)};\dots; D^{(n)}}
\left(
\begin{array}{c}
[(a):\theta^{(1)},\dots,\theta^{(n)}]
 : [(b^{{1}}):\phi^{(1)}];\dots;[(b^{{n}}):\phi^{(n)} ]
   \\
\, [(c):\psi^{(1)},\dots,\psi^{(n)}]
  : [(d^{{1}}):\delta^{(1)}];\dots;[(d^{{n}}):\delta^{(n)} ]
\end{array}
\, q; x_1,\dots, x_n \right)
\nonumber\\
&=&\sum_{s_1,\dots,s_n=0}^{\infty}\Omega(s_1,\dots,s_n)
\frac{x_1^{s_1}}{(q,q)_{s_1}}\dots \frac{x_n^{s_n}}{(q,q)_{s_n}}
\, ,
\label{Eq:Sriv:Def}
\end{eqnarray}
where
\begin{eqnarray}
\Omega(s_1,\dots,s_n)
=\frac{\prod_{j=1}^A (a_j,q)_{s_1\theta_j^{(1)}+\dots+ s_n\theta_j^{(n)}}  \prod_{j=1}^{B^{(1)}} (b^{(1)}_j,q)_{s_1\phi_j^{(1)}}\dots \prod_{j=1}^{B^{(n)}} (b^{(n)}_j,q)_{s_n\phi^{(n)}_j} }
      {\prod_{j=1}^C (c_j,q)_{s_1\psi_j^{(1)}+\dots+ s_n\psi_j^{(n)}}  \prod_{j=1}^{D^{(1)}} (d^{(1)}_j,q)_{s_1\delta_j^{(1)}}\dots \prod_{j=1}^{D^{(n)}} (d^{(n)}_j,q)_{s_n\delta^{(n)}_j} }
\, ,
\label{eq:def:omeg}
\end{eqnarray}
and all parameters
\begin{gather}
\label{Eq:Par}
\theta_1^{{(1)}},\psi_1^{{(1)}},\phi_1^{{(1)}},\delta_1^{{(1)}},\dots, \theta_A^{{(n)}},\psi_C^{{(n)}},\phi_{B^{(n)}}^{{(n)}},\delta_{D^{(n)}}^{{(n)}}
\end{gather}
are so constrained that multiple series (\ref{Eq:Sriv:Def}) converges.
Here we impose additional condition that parameters (\ref{Eq:Par}) are non-negative real numbers.

Here $(a,q)_n$ is q-shifted factorial and is defined by \cite{book1}:
\begin{gather}
(a,q)_n=\Pi_{m=0}^{n-1}(1-a q^m), \quad (a,q)_0=1,
\\
(a,q)_n=\frac{(a,q)_\infty}{(aq^n,q)_\infty}
\end{gather}
The last definition of q-shifted factorial helps us to extend the  definition area of variable  $n$ to negative integers and real values:
\begin{gather}
(a,q)_{-n}=\frac{1}{\Pi_{m=1}^{n}(1-a q^{-m})}=\frac{1}{(aq^{-n},q)_n}=\frac{(-q/a)^n q^{
%
%
\frac{n(n-1)}{2}
    }}{(q/a,q)_n}
\end{gather}
The q-analog of n, also known as the q-bracket or q-number of n, is defined to be
\begin{gather}
[n]_q=1+q^2+\dots+q^{n-1},
\end{gather}
The equivalent definition as
\begin{gather}
[n]_q=\frac{1-q^n}{1-q},
\end{gather}
gives us possibility to extend the definition over real values of $n$.

The  Jackson \cite{Jackson} q-derivative is a q-analog of the ordinary derivative of a function $f(x)$  defined as:
\begin{gather}
D_{x,q}f(x)=\frac{f(qx)-f(x)}{(q-1)x}.
\label{eq:derRule}
\end{gather}
It has product rule analogous to the ordinary derivative product, we could produce from above equation  the q-derivative for the ratio of two functions:
\begin{gather}
D_{x,q}=\frac{f(x)}{g(x)}=\frac{g(x)D_{x,q}f(x)-f(x)D_{x,q}g(x)}{g(qx)g(x)},
\end{gather}
more properties could be find in \cite{Diff1},\cite{book1},\cite{Diff3}.

For our purposes we shorten notation for  (\ref{Eq:Sriv:Def},\ref{eq::defH3}) functions to  names without parameters and variables and   explicitly written out only
those parameters that are different from those in definition (\ref{Eq:Sriv:Def},\ref{eq::defH3}). For example, $F(\lambda x_2)$ means
that instead of variable $x_2$ we will use variable $\lambda x_2$ in definition (\ref{Eq:Sriv:Def}) of $F$.


\section{Q-difference of q-analog of non-confluent Horn type hypergeometric function  $H_3$}
\label{sec3}
The q-analog of non-confluent Horn-type hypergeometric function $H_3(a,b,c,z_1,z_2)$
\begin{gather}
H_3(a,b,c,q,z_1,z_2)=\sum_{n_1,n_2}\frac{(a)_{2n_1+n2}(b)_{n_2}}{n_1!n_2! (c)_{n_1+n_2}}z_1^{n_1}z_2^{n_2}
\end{gather}
could be written in the form:
\begin{gather}
H_{q,3}(a,b,c,q,z_1,z_2)=\sum_{n_1,n_2}\frac{(a,q)_{2n_1+n_2}(b,q)_{n_2}}{(q,q)_{n_1}(q,q)_{n_2} (c,q)_{n_1+n_2}}z_1^{n_1}z_2^{n_2}.
\label{eq::defH3}
\end{gather}

\subsection{Q-derivative with respect to the upper parameters}

Derivative over upper one-summation parameter $b$ could be calculated by the same algorithm as in \cite{Ghany}:
\begin{gather}
D_{b,q}H_{q,3}(a,b,c,q,z_1,z_2)=\sum_{n_1,n_2}
\frac{(bq,q)_{n_2}-(b,q)_{n_2}}{(q-1)b }
\frac{(a,q)_{2n_1+n_2}}{(q,q)_{n_1}(q,q)_{n_2} (c,q)_{n_1+n_2}}z_1^{n_1}z_2^{n_2}
\nonumber\\
=\sum_{n_1,n_2}\frac{1}{(q-1)b }
\left(\frac{(bq,q)_\infty}{(bq^{n_2+1},q)_\infty}-\frac{(b,q)_\infty}{(bq^{n_2},q)_\infty}    \right)
\frac{(a,q)_{2n_1+n_2}}{(q,q)_{n_1}(q,q)_{n_2} (c,q)_{n_1+n_2}}z_1^{n_1}z_2^{n_2}
\nonumber\\
=\sum_{n_1,n_2}\frac{1}{(q-1)b }\frac{(b,q)_\infty}{(bq^{n_2},q)_\infty}
\left(\frac{1-bq^{n_2}}{1-b}-1    \right)
\frac{(a,q)_{2n_1+n2}}{(q,q)_{n_1}(q,q)_{n_2} (c,q)_{n_1+n_2}}z_1^{n_1}z_2^{n_2}
\nonumber\\
=\sum_{n_1,n_2}\frac{-1}{(1-b) }\frac{q^{n_2}-1}{(q-1)}
\frac{(a,q)_{2n_1+n_2}(b,q)_{n_2}}{(q,q)_{n_1}(q,q)_{n_2} (c,q)_{n_1+n_2}}z_1^{n_1}z_2^{n_2}.
\label{eq:DerH3:b}
\end{gather}
By using the derivative of $H_3$ over variable $z_i$
%
\begin{gather}
D_{z_i,q}H_{q,3}(a,b,c,q,z_1,z_2)=
\sum_{n_1,n_2}
\frac{1}{z_i}[n_i]_q
\frac{(a,q)_{2n_1+n_2}(b,q)_{n_2}}{(q,q)_{n_1}(q,q)_{n_2} (c,q)_{n_1+n_2}}z_1^{n_1}z_2^{n_2},
\label{eq:DerSimple}
\end{gather}
we could write (\ref{eq:DerH3:b}) as follows:
\begin{gather}
D_{b,q}H_{q,3}(a,b,c,q,z_1,z_2)=-\frac{z_2}{1-b}D_{z_2,q}H_{q,3}(a,b,c,q,z_1,z_2).
\label{H3::Der::Up}
\end{gather}

It is not so straightforward calculation of  q-derivative with respect to the double-summation parameter $a$ as in the case of one-summation index $b$, eq. (\ref{eq:DerH3:b}):
\begin{gather}
D_{a,q}H_{q,3}(a,b,c,q,z_1,z_2)=\sum_{n_1,n_2}\frac{(aq,q)_{2n_1+n_2}-(a,q)_{2n_1+n_2}}{(q-1)a }\frac{(b,q)_{n_2}}{(q,q)_{n_1}(q,q)_{n_2} (c,q)_{n_1+n_2}}z_1^{n_1}z_2^{n_2}
\nonumber\\
=\sum_{n_1,n_2}\frac{1}{(q-1)a }
\left(\frac{(aq,q)_\infty}{(aq^{2n_1+n_2+1},q)_\infty}-\frac{(a,q)_\infty}{(aq^{2n_1+n_2},q)_\infty}    \right)
\frac{(b,q)_{n2}}{(q,q)_{n_1}(q,q)_{n_2} (c,q)_{n_1+n_2}}z_1^{n_1}z_2^{n_2}
\nonumber\\
=\sum_{n_1,n_2}\frac{1}{(q-1)a }\frac{(a,q)_\infty}{(aq^{2n_1+n_2},q)_\infty}
\left(\frac{1-aq^{2n_1+n_2}}{1-a}-1    \right)
\frac{(b,q)_{n2}}{(q,q)_{n_1}(q,q)_{n_2} (c,q)_{n_1+n_2}}z_1^{n_1}z_2^{n_2}
\nonumber\\
=\sum_{n_1,n_2}\frac{-1}{(1-a) }\frac{q^{2n_1+n_2}-1}{(q-1)}
\frac{(a,q)_{2n_1+n2}(b,q)_{n_2}}{(q,q)_{n_1}(q,q)_{n_2} (c,q)_{n_1+n_2}}z_1^{n_1}z_2^{n_2}.
\label{eq:DerH3:a}
\end{gather}
Unlike to the eq. (\ref{eq:DerH3:b})  we have here  q-number $[2n_1+n_2]_q$, which contains double summation indexes with the integer coefficient 2.
By means of identity that splits double summation indexes to two  items with t one summation indexes
\begin{gather}
[2n_1+n_2]_q=\frac{1}{2}\left(
(1+q^{2n_1})[n_2]_q+(1+q^{n_2})[2n_1]
\right),
\label{eq:sumSplit}
\end{gather}
and the expression for the q-derivative of $H_3$ function with squared  variable $z_2$,
\begin{gather}
D_{z_2,q}H_{q,3}(a,b,c,q,z_1, z_2^2)=\sum_{n_1,n_2}\frac{1}{z_2}[2n_2]_q\frac{(a,q)_{2n_1+n_2}(b,q)_{n_2}}{(q,q)_{n_1}(q,q)_{n_2} (c,q)_{n_1+n_2}}z_1^{n_1}z_2^{2 n_2},
\end{gather}
we get the following answer for q-derivative with respect to the double-summation index parameter $a$:
\begin{eqnarray}
D_{a,q}H_{q,3}(a,b,c,q,z_1,z_2)
\!\!\!\!\!
&=&\frac{-1}{2(1-a)}
\nonumber
\\
&\times&
\!\!\!\!\!
\left[
z_2 D_{z_2,q}(H_{q,3}+H_{q,3}(q^2 z_1))
+z_1 D_{z_1,q}(H_{q,3}(z_1^2)+H_{q,3}(z_1^2,q z_2))
\right]
\!
.
\end{eqnarray}

\subsection{Q-derivative with respect to the lower parameters}

The similar procedure could be applied for derivation with respect to the double-summation lower parameter $c$:
\begin{eqnarray}
D_{c,q}H_{q,3}(a,b,c,z_1,z_2)&&=\sum_{n_1,n_2}
\frac{-1}{(c,q)_{n_1+n_2}}
\frac{(cq,q)_{n_1+n_2}-(c,q)_{n_1+n_2}}{(q-1)c }
\nonumber\\
&&\times
\frac{(a,q)_{2n_1+n_2}(b,q)_{n_2}}{(q,q)_{n_1}(q,q)_{n_2}(cq,q)_{n_1+n_2} }z_1^{n_1}z_2^{n_2}
\nonumber\\
&&=\sum_{n_1,n_2}
\frac{-1}{c(q-1)(c,q)_{n_1+n_2}}
\left(\frac{(cq,q)_\infty}{(cq^{n_1+n_2+1},q)_\infty}-\frac{(c,q)_\infty}{(cq^{n_1+n_2},q)_\infty}    \right)
\nonumber\\
&&\times
\frac{(a,q)_{2n_1+n_2}(b,q)_{n2}}{(q,q)_{n_1}(q,q)_{n_2} (cq,q)_{n_1+n_2}}z_1^{n_1}z_2^{n_2}
\nonumber\\
&&=\sum_{n_1,n_2}\frac{-1}{(q-1)c }
\left(\frac{1-cq^{n_1+n_2}}{1-c}-1    \right)
\nonumber\\
&&\times
\frac{(a,q)_{2n_1+n_2}(b,q)_{n2}}{(q,q)_{n_1}(q,q)_{n_2} (cq,q)_{n_1+n_2}}z_1^{n_1}z_2^{n_2}
\nonumber\\
&&=\sum_{n_1,n_2}\frac{1}{(1-c) }\frac{q^{n_1+n_2}-1}{(q-1)}
\frac{(a,q)_{2n_1+n2}(b,q)_{n_2}}{(q,q)_{n_1}(q,q)_{n_2} (cq,q)_{n_1+n_2}}z_1^{n_1}z_2^{n_2}.
\label{eq:DerH3:Lower}
\end{eqnarray}
With the help of identity similar to the eq. (\ref{eq:sumSplit})
\begin{gather}
[n_1+n_2]_q=\frac{1}{2}\left(
(1+q^{n_1})[n_2]_q+(1+q^{n_2})[n_1]
\right),
\label{eq:sumSplitSimple}
\end{gather}
 and eq. (\ref{eq:DerSimple}) we could write the final answer for the q-derivative with the respect to the lower parameter  $c$:
\begin{gather}
D_{c,q}H_{q,3}(a,b,c,z_1,z_2)=\frac{1}{2(1-c)}\left[
z_2 D_{z_2,q}(H_{q,3}+H_{q,3}(q z_1))
+z_1 D_{z_1,q}(H_{q,3}+H_{q,3}(q z_2))
\right].
\label{H3::Der::Down}
\end{gather}

\section{Q-derivative of q-extension of the generalized Lauricella series  with respect to the parameters}
\label{sec:Main}
\subsection{Derivative with respect to the upper parameter }
For calculation of the derivative we need and extension of identity (\ref{eq:sumSplit}) to the case of multiple summation index with real non-negative coefficients $\theta_i$.
We will put it here explicitly:
\begin{gather}
[\theta_1m_1+\dots\theta_k m_k]_q=\frac{1}{k}( (1+q^{\theta_2 m_2}+q^{\theta_2 m_2+\theta_3 m_3}+...+q^{\theta_2 m_2+\dots\theta_k m_k})[\theta_1 m_1]_q
\nonumber
\\
+(1+q^{\theta_3 m_3}+q^{\theta_3 m_3+\theta_4 m_4}+...+q^{\theta_3 m_3+\dots\theta_k m_k+\theta_1m_1})[\theta_2 m_2]_q
\nonumber
\\
\cdots
\nonumber
\\
+(1+q^{\theta_1 m_1}+q^{\theta_1 m_1+\theta_2 m_2}+...+q^{\theta_1 m_1+\dots\theta_{k-1} m_{k-1}})[\theta_k m_k]_q
),
\label{eq:genSum}
\end{gather}
here all  $m_j$ are non-negative integers.

The derivative over variable $z_k$ of hypergeometric function (\ref{Eq:Sriv:Def}) reads:
%
\begin{gather}
D_{z_k,q}F(z_1,...,z_n)=
\sum_{s_1,\dots,s_n=0}^{\infty}
\frac{1}{z_k}[s_k]_q
\Omega(s_1,\dots,s_n)
\frac{z_1^{s_1}}{(q,q)_{s_1}}\dots \frac{z_n^{s_n}}{(q,q)_{s_n}}.
\label{eq:genDerX}
\end{gather}

Now we could obtain the explicit formula for q-derivative of hypergeometric function $F$ with respect to the upper parameter $a_j$:
\begin{gather}
D_{a_j,q}F=\sum_{s_1,\dots,s_n=0}^{\infty}
\frac{(a_jq,q)_{s_1\theta_j^{(1)}+\dots+ s_n\theta_j^{(n)}}-(a_j,q)_{s_1\theta_j^{(1)}+\dots+ s_n\theta_j^{(n)}}}{(q-1)a_j }
\nonumber
\\
\times\frac{1}{(a_j,q)_{s_1\theta_j^{(1)}+\dots+ s_n\theta_j^{(n)}}}
\Omega(s_1,\dots,s_n)
\frac{x_1^{s_1}}{(q,q)_{s_1}}\dots \frac{x_n^{s_n}}{(q,q)_{s_n}}
\nonumber
\\
=\sum_{s_1,\dots,s_n=0}^{\infty}
\frac{1}{(q-1)a_j }\frac{(a_j,q)_\infty}{(a_jq^{s_1\theta_j^{(1)}+\dots+ s_n\theta_j^{(n)}},q)_\infty}
\left(\frac{1-a_jq^{s_1\theta_j^{(1)}+\dots+ s_n\theta_j^{(n)}}}{1-a_j}-1    \right)
\nonumber
\\
\times\frac{1}{(a_j,q)_{s_1\theta_j^{(1)}+\dots+ s_n\theta_j^{(n)}}}
\Omega(s_1,\dots,s_n)
\frac{x_1^{s_1}}{(q,q)_{s_1}}\dots \frac{x_n^{s_n}}{(q,q)_{s_n}}
\nonumber
\\
=\sum_{s_1,\dots,s_n=0}^{\infty}
\frac{-1}{(1-a_j) }\frac{q^{s_1\theta_j^{(1)}+\dots+ s_n\theta_j^{(n)}}-1}{(q-1)}
\nonumber
\\
\times
\Omega(s_1,\dots,s_n)
\frac{x_1^{s_1}}{(q,q)_{s_1}}\dots \frac{x_n^{s_n}}{(q,q)_{s_n}},
\end{gather}
by applying the multivariable versions (\ref{eq:genSum}, \ref{eq:genDerX})   of (\ref{eq:sumSplit}, \ref{eq:DerSimple})
we obtain the main result for q-derivation of   the generalized Lauricella q-extension series with respect to the upper parameter:
\begin{gather}
D_{a_j,q}F=\frac{-1}{k(1-a_j)}\left[
 z_1 D_{z_1,q}(F(z_1^{\theta_1})+ F(z_1^{\theta_1},q^{\theta_2}z_2)+ F(z_1^{\theta_1},q^{\theta_2}z_2,q^{\theta_3}z_3)
 \right.
 \nonumber
 \\
 \left.
 +\dots+F(z_1^{\theta_1},q^{\theta_2}z_2,\dots,q^{\theta_k}z_k))
 \right.
\nonumber
\\
\left.
+z_2D_{z_2,q}(F(z_2^{\theta_2})+F(z_2^{\theta_2},q^{\theta_3}z_3)+F(z_2^{\theta_2},q^{\theta_3}z_3,q^{\theta_4}z_4)
\right.
\nonumber
\\
\left.
+\dots
+F(q^{\theta_1}z_1,z_2^{\theta_2},q^{\theta_3}z_3,\dots,q^{\theta_k}z_k))
\right.
\nonumber
\\
\cdots
\nonumber
\\
\left.
+z_kD_{z_k,q}(F(z_k^{\theta_k})+F(q^{\theta_1}z_1, z_k^{\theta_k})+F(q^{\theta_1}z_1,q^{\theta_2}z_2, z_k^{\theta_k})
\right.
\nonumber
\\
\left.
+\dots+F(q^{\theta_1}z_1,\dots, q^{\theta_{k-1}}z_{k-1}, z_k^{\theta_k}))
\right].
\label{eq:DerA}
\end{gather}

For the one-summation upper  parameter index q-derivation $b_j^{(1)}=b_j$ we could produce:
\begin{gather}
D_{b_j,q}F=\sum_{s_1,\dots,s_n=0}^{\infty}
\frac{(b_jq,q)_{s_1\phi_j^{(1)}}-(b_j,q)_{s_1\phi_j^{(1)}}}{(q-1)b_j }
\frac{1}{(b_j,q)_{s_1\phi_j^{(1)}}}
\Omega(s_1,\dots,s_n)
\frac{x_1^{s_1}}{(q,q)_{s_1}}\dots \frac{x_n^{s_n}}{(q,q)_{s_n}}
\nonumber
\\
=\sum_{s_1,\dots,s_n=0}^{\infty}
\frac{1}{(q-1)b_j }\frac{(b_j,q)_\infty}{(b_jq^{s_1\phi_j^{(1)}},q)_\infty}
\left(\frac{1-b_jq^{s_1\phi_j^{(1)}}}{1-b_j}-1    \right)
\nonumber
\\
\times
\frac{1}{(b_j,q)_{s_1\phi_j^{(1)}}}
\Omega(s_1,\dots,s_n)
\frac{x_1^{s_1}}{(q,q)_{s_1}}\dots \frac{x_n^{s_n}}{(q,q)_{s_n}}
\nonumber
\\
=\sum_{s_1,\dots,s_n=0}^{\infty}
\frac{-1}{(1-b_j) }\frac{q^{s_1\phi_j^{(1)}}-1}{(q-1)}
\Omega(s_1,\dots,s_n)
\frac{x_1^{s_1}}{(q,q)_{s_1}}\dots \frac{x_n^{s_n}}{(q,q)_{s_n}},
\end{gather}
and with the help of eq. (\ref{eq:genDerX}) we obtain:
\begin{gather}
D_{b_j,q}F=\frac{-1}{(1-b_j)}
 z_1 D_{z_1,q}F(q^{\phi_j^{(1)}}z_1).
\label{eq:DerB}
\end{gather}
The eq. (\ref{eq:DerB}) also could be obtained by applying $a_j\to b_j$, $\theta_j^{(1)}\to\phi_j^{(1)}$, $\theta_j^{(2)}\dots\theta_j^{(n)}\to 0$.

The derivative with respect to other upper parameters $b_j^{(2)}\dots b_j^{(n)}$ could be obtained from eq. ({\ref{eq:DerB}}) by substitution
$b_j\to b_j^{(i)}$, $\phi_j^{(1)}\to \phi_j^{(i)}$, $z_1\to z_i$:
\begin{gather}
D_{b_j^{(i)},q}F=\frac{-1}{(1-b_j^{(i)})} z_i D_{z_i,q}F(q^{\phi_j^{(i)}}z_i).
\end{gather}

\subsection{Q-derivative with respect to the lower parameter }

With the same procedure we consider the q-derivative of the Srivastava and Daoust multivariable hypergeometric  series with respect to the lower parameter $c_j$:
\begin{gather}
D_{c_j,q}F=\sum_{s_1,\dots,s_n=0}^{\infty}
\frac{(c_jq,q)_{s_1\psi_j^{(1)}+\dots+ s_n\psi_j^{(n)}}-(a_j,q)_{s_1\psi_j^{(1)}+\dots+ s_n\psi_j^{(n)}}}{(q-1)c_j }
\nonumber
\\
\times\frac{-1}{(qc_j,q)_{s_1\psi_j^{(1)}+\dots+ s_n\psi_j^{(n)}}}
\Omega(s_1,\dots,s_n)
\frac{x_1^{s_1}}{(q,q)_{s_1}}\dots \frac{x_n^{s_n}}{(q,q)_{s_n}}
\nonumber
\\
=\sum_{s_1,\dots,s_n=0}^{\infty}
\frac{-1}{(q-1)c_j }\frac{(c_j,q)_\infty}{(c_jq^{s_1\psi_j^{(1)}+\dots+ s_n\psi_j^{(n)}},q)_\infty}
\left(\frac{1-c_jq^{s_1\psi_j^{(1)}+\dots+ s_n\psi_j^{(n)}}}{1-a_j}-1    \right)
\nonumber
\\
\times\frac{1}{(qc_j,q)_{s_1\psi_j^{(1)}+\dots+ s_n\psi_j^{(n)}}}
\Omega(s_1,\dots,s_n)
\frac{x_1^{s_1}}{(q,q)_{s_1}}\dots \frac{x_n^{s_n}}{(q,q)_{s_n}}
\nonumber
\\
=\sum_{s_1,\dots,s_n=0}^{\infty}
\frac{1}{(1-c_j) }\frac{q^{s_1\psi_j^{(1)}+\dots+ s_n\psi_j^{(n)}}-1}{(q-1)}
\frac{(c_j,q)_{s_1\psi_j^{(1)}+\dots+ s_n\psi_j^{(n)}}}{(qc_j,q)_{s_1\psi_j^{(1)}+\dots+ s_n\psi_j^{(n)}}}
\nonumber
\\
\times
\Omega(s_1,\dots,s_n)
\frac{x_1^{s_1}}{(q,q)_{s_1}}\dots \frac{x_n^{s_n}}{(q,q)_{s_n}}.
\end{gather}
By applying  (\ref{eq:genSum}, \ref{eq:genDerX})   we obtain:
\begin{gather}
D_{c_j,q}F=\frac{1}{k(1-c_j)}\left[
 z_1 D_{z_1,q}(F(qc_j,z_1^{\theta_1})+ F(qc_j,z_1^{\theta_1},q^{\theta_2}z_2)+ F(qc_j,z_1^{\theta_1},q^{\theta_2}z_2,q^{\theta_3}z_3)
 \right.
 \nonumber
 \\
 \left.
 +\dots+F(qc_j,z_1^{\theta_1},q^{\theta_2}z_2,\dots,q^{\theta_k}z_k))
 \right.
\nonumber
\\
\left.
+z_2D_{z_2,q}(F(qc_j,z_2^{\theta_2})+F(qc_j,z_2^{\theta_2},q^{\theta_3}z_3)+F(qc_j,z_2^{\theta_2},q^{\theta_3}z_3,q^{\theta_4}z_4)
\right.
\nonumber
\\
\left.
+\dots
+F(qc_j,q^{\theta_1}z_1,z_2^{\theta_2},q^{\theta_3}z_3,\dots,q^{\theta_k}z_k))
\right.
\nonumber
\\
\cdots
\nonumber
\\
\left.
+z_kD_{z_k,q}(F(qc_j,z_k^{\theta_k})+F(qc_j,q^{\theta_1}z_1, z_k^{\theta_k})+F(qc_j,q^{\theta_1}z_1,q^{\theta_2}z_2, z_k^{\theta_k})
\right.
\nonumber
\\
\left.
+\dots+F(qc_j,q^{\theta_1}z_1,\dots, q^{\theta_{k-1}}z_{k-1}, z_k^{\theta_k}))
\right].
\label{eq:DerC}
\end{gather}

For the one-summation index lower  parameter  q-derivation $d_j^{(i)}$ we could produce:
\begin{gather}
D_{d_j^{(i)},q}F=\sum_{s_1,\dots,s_n=0}^{\infty}
\frac{(d_j^{(i)}q,q)_{s_i\delta_j^{(i)}}-(d_j^{(i)},q)_{s_i\delta_j^{(i)}}}{(q-1)d_j^{(i)} }
\frac{-1}{(d_j^{(i)}q,q)_{s_i\delta_j^{(i)}}}
\Omega(s_1,\dots,s_n)
\frac{x_1^{s_1}}{(q,q)_{s_1}}\dots \frac{x_n^{s_n}}{(q,q)_{s_n}}
\nonumber
\\
=\sum_{s_1,\dots,s_n=0}^{\infty}
\frac{1}{(q-1)d_j^{(i)} }\frac{(d_j^{(i)},q)_\infty}{(d_j^{(i)}q^{s_i\delta_j^{(i)}},q)_\infty}
\left(\frac{1-d_j^{(i)}q^{s_i\delta_j^{(i)}}}{1-d_j^{(i)}}-1    \right)
\nonumber
\\
\times
\frac{-1}{(d_j^{(i)}q,q)_{s_i\delta_j^{(i)}}}
\Omega(s_1,\dots,s_n)
\frac{x_1^{s_1}}{(q,q)_{s_1}}\dots \frac{x_n^{s_n}}{(q,q)_{s_n}}
\nonumber
\\
=\sum_{s_1,\dots,s_n=0}^{\infty}
\frac{1}{(1-d_j^{(i)}) }\frac{q^{s_i\delta_j^{(i)}}-1}{(q-1)}
\frac{(d_j^{(i)},q)_{s_i\delta_j^{(i)}}}{(d_j^{(i)}q,q)_{s_i\delta_j^{(i)}}}
\Omega(s_1,\dots,s_n)
\frac{x_1^{s_1}}{(q,q)_{s_1}}\dots \frac{x_n^{s_n}}{(q,q)_{s_n}},
\end{gather}
and with the help of eq. (\ref{eq:genDerX}) we obtain:
\begin{gather}
D_{d_j^{(i)},q}F=\frac{1}{(1-d_j^{(i)})}
 z_i D_{z_i,q}F(q^{\delta_j^{(i)}}z_i).
\label{eq:DerD}
\end{gather}

The equations (\ref{eq:DerA}, \ref{eq:DerB}, \ref{eq:DerC}, \ref{eq:DerD}) show the final expressions for any q-derivative with respect to the parameters of the   q-extension of the generalized Lauricella series.

\section{Conclusion}
We extend the algorithm \cite{Bytev:2017jmx} for derivation of generalized  Lauricella hypergeometric series  to the case of basic hypergeometric function q-derivatives and consider the general case of q-derivatives with respect to the parameters for q-extension of the generalized Lauricella series. Namely, the  derivation formula for upper parameter is presented in (\ref{eq:DerA},\ref{eq:DerB}), derivatives over lower parameter could be find at (\ref{eq:DerC}, \ref{eq:DerD}). By that equations we cover all possible cases  parameter derivation of Srivastava and Daoust basic multivariable hypergeometric function with positive summation index  and explicitly show that the parameter derivative could be expressed as a finite sum of derivative over variables with  shifted values.

As an example we produce q-derivative of   non-confluent Horn-type hypergeometric function  $H_3(a,b,c,z_1,z_2)$  q-analog, over upper (\ref{H3::Der::Up}) and lower  (\ref{H3::Der::Down}) parameters correspondingly.

\section*{Acknowledgments}

This work was partially supported by the National Key Research and Development Program of China (No. 2016YFE0130800), the National Natural Science Foundation of China (Grants No. 11975320) and by the Heisenberg--Landau Program.



\begin{thebibliography}{99}


\bibitem{Sriv1}
H.~M.~Srivastava,
\newblock {\it Certain q-polynomial expansions for functions of several variables},
\newblock IMA \ J.\ App. \ Math. {\bf 30} 315 (1983).

\bibitem{book1}
G.~Gasper, M.~Rahman
\newblock {\it Basic Hypergeometric Series},
\newblock  Cambridge, Encyclopedia of Mathematics and its Applications) (2004).

\bibitem{Ancarani0}
L.~U. Ancarani and G.~Gasaneo,
\newblock {\it Derivatives of any order of the Gaussian
  hypergeometric function $_2F_1(a,b,c;z)$ with respect to the parameters $a$, $b$ and $c$},
\newblock Journal of Physics A: Mathematical and Theoretical {\bf 42}, 395208 (2009).

\bibitem{Ancarani1}
L.~U. Ancarani and G.~Gasaneo,
\newblock {\it Derivatives of any order of the hypergeometric function $_pF_q(a_1, ..., a_p; b_1, ..., b_q;z)$ with
  respect to the parameters $a_i$ and $b_i$},
\newblock Journal of Physics A: Mathematical and Theoretical {\bf 43}, 085210 (2010).

\bibitem{Ancarani2}
L.~U. Ancarani, J.~A.~D. Punta, and G.~Gasaneo,
\newblock {\it Derivatives of Horn hypergeometric functions with respect to their parameters},
\newblock J.\ Math.\ Phys.\ {\bf 58}, 073504 (2017).

\bibitem{Sahai}
V.~Sahai and A.~Verma,
\newblock {\it Derivatives of Appell functions with respect to parameters},
\newblock Journal of Inequalities and Special Functions {\bf 6}, 1 (2015).

\bibitem{Bujar}
B.~Xh.~Fejzullahu,
\newblock {\it Parameter derivatives of the generalized hypergeometric function},
\newblock Integral Transforms and Special Functions {\bf 28}, 781  (2017).

\bibitem{Adams:2015ydq}
L.~Adams, C.~Bogner and S.~Weinzierl,
J. Math. Phys. \textbf{57} (2016) no.3, 032304
doi:10.1063/1.4944722
[arXiv:1512.05630 [hep-ph]].

\bibitem{Adams:2016xah}
L.~Adams, C.~Bogner, A.~Schweitzer and S.~Weinzierl,
J. Math. Phys. \textbf{57} (2016) no.12, 122302
doi:10.1063/1.4969060
[arXiv:1607.01571 [hep-ph]].

\bibitem{Passarino:2016zcd}
G.~Passarino,
Eur. Phys. J. C \textbf{77} (2017) no.2, 77
doi:10.1140/epjc/s10052-017-4623-1
[arXiv:1610.06207 [math-ph]].


\bibitem{levin}
A.~LEVIN,  
Compositio Mathematica, \textbf{106}  (1997) no.3, 267
doi:10.1023/A:1000193320513


\bibitem{firstDer1}
M.~Abramowitz,
\newblock {\it Handbook of mathematical functions, with formulas, graphs, and mathematical tables},
\newblock (Dover Publications Inc., New York, N.Y., 1974).

\bibitem{firstDer2}
Yu.~A. Brychkov,
\newblock {\it Handbook of special functions: derivatives, integrals, series and other formulas},
\newblock (Chapman and Hall/CRC Press, New York, 2008).

\bibitem{firstDer3}
J.~Froehlich,
\newblock {\it Parameter derivatives of the Jacoby polynomials and the Gaussian hypergeometric function},
\newblock Integral Transforms and Special Functions {\bf 2}, 253 (1994).


\bibitem{Ghany}
H.~Ghany,
\newblock {\it Q-derivative of basic hypergeometric series with respect to parameters},
\newblock  Int. \ J. of\ Math.\ Analysis, {\bf 3} 1617 (2009).


\bibitem{Jackson}
F.~N.~Jackson,
\newblock {\it On q-functions and certain difference operator},
\newblock  Trans. \ Roy. \ Soc.\ Edin, {\bf 46}  (1908).

\bibitem{Diff1}
H.~Exton,
\newblock {\it Q-hypergeometric functions and applications},
\newblock  Ellis Horwood Limited, Chichester,  (1983).

\bibitem{Diff3}
F.~N.~Jackson,
\newblock {\it On basic hypergeometric functions},
\newblock  Quart.\ J.\ Math.\ {\bf 13} 69 (1942).


\bibitem{Bytev:2017jmx}
V.~V.~Bytev and B.~A.~Kniehl,
\newblock {\it Derivatives of Horn-type hypergeometric functions with respect to their parameters},
\newblock arXiv:1712.07579 [math-ph].
\newblock Nucl.\ Phys.\ B {\bf 952} (2020)
  
\bibitem{qref1}
H.~A.~Ghany,
\newblock {\it Q-derivative of basic hypergoemtric series with respect to parameters},
\newblock  Int.\ J.\ Math.\ Anal {\bf 3} 33 (2009).

\bibitem{qref2}
V.~B.~Kuznetsov and ~E.~K.~Sklyanin, 
\newblock {\it Factorisation of Macdonald polynomials, Symmetries and integrability of difference equations}
\newblock  London \ Math. \ Soc. \ Lecture \ Note \ Ser., {\bf 255 } 370 (1999).

\bibitem{qref3}
V.~Sahai and A.~Verma, 
\newblock {\it nth-order q-derivatives of multivariable q-hypergeometric series with respect to parameters},
\newblock   Asian-Eur.\ J.\ Math. {\bf 7}  29 (2014).

\bibitem{qref4}
V.~Sahai and A.~Verma, 
\newblock {\it Derivatives of Appell functions with respect to parameters},
\newblock J. \ Inequal.\ Spec.\ Funct., {\bf 6} 1 (2015)

\bibitem{qref5}
V.~Sahai and A.~Verma,
\newblock {\it Derivatives of Appell functions with respect to parameters},
\newblock KYUNGPOOK \ Math. \ J.  {\bf 56 } 911 (2016)

\end{thebibliography}
\end{document}